\documentclass[12pt, final]{amsart}
\usepackage[mathscr]{eucal}

\theoremstyle{plain}
\newtheorem*{theorem}{Theorem}
\newtheorem*{lemma}{Lemma}

\theoremstyle{remark}
\newtheorem*{remark}{Remark}

\begin{document}
\title[ON LOCAL AUTOMORPHISMS OF GROUP ALGEBRAS]
{ON LOCAL AUTOMORPHISMS OF GROUP ALGEBRAS OF COMPACT GROUPS}
\author[LAJOS MOLN\'AR]{LAJOS MOLN\' AR${}^1$}
\address{Institute of Mathematics\\
         Lajos Kossuth University\\
         4010 Debrecen, P.O.Box 12, Hungary}
\email{molnarl@math.klte.hu}
\author[BORUT ZALAR]{BORUT ZALAR${}^2$}
\address{Department of Basic Sciences\\
         Faculty of Civil Engineering\\
         University of Maribor, Smetanova 17\\
         2000 Maribor, Slovenia}
\email{borut.zalar@uni-mb.si}
\thanks{  This research was supported from the following sources:\\
          a)${}^{1,2}$ Joint Hungarian-Slovene research project
             supported
             by OMFB in Hungary and the Ministry of Science and
             Technology in Slovenia, Reg. No. SLO-2/96,\\
          b)${}^1$ Hungarian National Foundation for Scientific Research
             (OTKA), Grant No. T--016846 F--019322,\\
          c)${}^1$ A grant from the Ministry of Education, Hungary, Reg.
            No. FKFP 0304/1997,
          d)${}^2$ A grant from the Ministry of Science and
            Technology, Slovenia}
\date{February 4, 1998}
\begin{abstract}
We show that with few exceptions
every local isometric automorphism of the group algebra $L^p(G)$
of a compact metric group $G$ is an isometric automorphism.
\end{abstract}
\maketitle

In the last decade considerable work has been done concerning
certain local maps of operator algebras.
The originators of this research are Kadison and
Larson. In \cite{Kad}, Kadison studied
local derivations on a von Neumann algebra $\mathcal R$.
A continuous linear map on $\mathcal R$
is called a local derivation if it agrees with some
derivation at each point (the derivations possibly differring from point
to point) in the algebra.
This investigation was motivated by the study of Hochschild cohomology
of operator algebras. It was proved in \cite{Kad}
that in the above setting, every local derivation is a
derivation. Independently, Larson and Sourour proved in \cite{LaSo} that
the same conclusion holds true for local derivations of the full
operator algebra $\mathscr B(\mathscr X)$, where $\mathscr X$ is a
Banach space. For other results on local derivations
of various algebras see, for example, \cite{Bre, BrSe1, Cri, Shu,
ZhXi}. Besides derivations,
there is at least one additional very important class of transformations
on Banach algebras which certainly deserves attention. This is the
group of automorphisms. In
\cite[Some concluding remarks (5), p. 298]{Lar}, from the
view-point of reflexivity, Larson raised the problem of local
automorphisms (the definition should be self-explanatory) of Banach
algebras. In his joint paper with Sourour \cite{LaSo},
it was proved that if $\mathscr X$ is an infinite dimensional Banach
space, then every surjective local automorphism of $\mathscr B(\mathscr
X)$ is an automorphism (see also \cite{BrSe1}).
For a separable infinite dimensional Hilbert
space $\mathscr H$, it was shown in \cite{BrSe2} that the above
conclusion
holds true without the assumption on surjectivity, i.e. every local
automorphism of $\mathscr B(\mathscr H)$ is an automorphism.
For other results on local automorphisms
of various operator algebras, we refer to \cite{BaMo, MolStud, MolLond,
MolJFA}.

In this note we investigate a similar problem for the
$L^p$-algebras (convolution algebras) of compact
groups which are of fundamental importance in harmonic analysis.
For the local isometric automorphisms of $L^p(G)$ (i.e. linear maps on
$L^p(G)$ which agree with some isometric automorphism at each point
in the algebra) we obtain the following result.

\begin{theorem}
Let $G$ be a first countable compact group
and let $1\leq p \leq \infty$. Then in the following cases
every local isometric automorphism of $L^p(G)$ is an
isometric automorphism.
\begin{itemize}
\item[(i)]  $p\neq 2$ and either $G$ is commutative,
            or $G$ is finite, or
            $G$ has an at least three-dimensional irreducible
            representation.
\item[(ii)] $p=2$ and $G$ does not have any two-dimensional irreducible
            representation.
\end{itemize}
\end{theorem}

\begin{remark}
Concerning the topological condition in our result we note that
by \cite[(8.3) Theorem]{HeRo} every first countable
compact group is metrizable.
\end{remark}

It is our approach to study first the restriction of our local isometric
automorphism $\psi$, acting on $L^p(G)$, to the subalgebra $C(G)$ of all
continuous complex valued functions on $G$. We equip $C(G)$ with the
usual supremum norm. In the proof of our theorem we use the following

\begin{lemma}
Let $X$ be a first countable compact Hausdorff space. If
$\psi:C(X) \to C(X)$ is a local surjective isometry, then
it is a surjective isometry.
\end{lemma}

\begin{proof}
By Banach-Stone theorem, every surjective isometry of
$C(X)$ is of the form $f\mapsto \tau \cdot f\circ \varphi$, where $\tau
:X \to \mathbb C$ is a continuous function of modulus 1 and $\varphi:
X\to X$ is a homeomorphism. It is now apparent that $\psi$ sends
continuous functions of modulus 1 to functions of the same kind.
Therefore, $\psi$ preserves the unitary elements of the $C^*$-algebra
$C(X)$. A result of Russo and Dye
\cite[Corollary 2]{RuDy} says that in that case $\psi$ is a Jordan
*-homomorphism followed by multiplication by a fixed unitary
element.
Without any loss of generality we may assume that $\psi(1)=1$.
Thus, we obtain that $\psi$ is an endomorphism of $C(X)$.
It is a folk result that every endomorphism of $C(X)$ which sends 1
to 1 is of the form
\[
f \longmapsto f\circ \varphi_0
\]
where $\varphi_0:X\to X$ is a continuous function. Since $\psi$ is an
isometry, it readily follows from Urysohn's lemma
that $\varphi_0$ is surjective. It remains to prove
that $\varphi_0$ is injective as well. To this end, suppose on the
contrary that there are different points $x,y \in X$ such that
$\varphi_0(x)=\varphi_0(y)=z$.
We construct a continuous function $f:X \to \mathbb C$ as follows.
Let $(U_n)$ be a monotone decreasing sequence of open sets in $X$ such
that $\cap_n U_n=\{ z\}$. By Urysohn's lemma, for every $n$ we have a
continuous function $f_n:X \to [0,1]$ such that
\[
f_n(z)=1 \enskip \text{and} \enskip f_n (t)=0 \enskip (t\in X\setminus
U_n).
\]
Let
\[
f=1-\sum_{n=1}^\infty  \frac{1}{2^n} f_n.
\]
Clearly, $f$ is a continuous function and
$f(t)=0$ if and only if $t=z$. Since $\psi$
is a local surjective isometry, there exist
a continuous function $\tau :X \to \mathbb C$ of modulus 1
and a homeomorphism $\varphi :X\to X$ such that
\[
f\circ \varphi_0 =\psi(f)=\tau \cdot f\circ \varphi.
\]
It follows that
\[
\{ x, y\} \subset (f\circ \varphi_0)^{-1}(0)=(f\circ
\varphi)^{-1}(0)=\varphi^{-1}(z).
\]
Since $\varphi$ is bijective, this is a contradiction. Consequently,
$\varphi_0$ is injective which implies the surjectivity of $\psi$.
\end{proof}

For the proof of our theorem we also need the following well-known
facts from harmonic analysis (see \cite[Sections 27,28,31]{HeRo}).

Let $\mathsf{X}$ be the dual object of $G$. For every $\sigma\in
\mathsf{X}$ let $d_\sigma$ be the dimension of the irreducible
representations of $G$ belonging to the equivalence class $\sigma$. Let
$\mathfrak{C}_0(\mathsf{X})$ denote the subset of
$\prod_{\sigma\in \mathsf{X}} M_{d_\sigma}(\mathbb C)$
consisting of those elements
$(E_\sigma)_{\sigma\in \mathsf{X}}$ which
vanish at infinity (here we consider the operator norm of matrices).
Clearly, $\mathfrak{C}_0(\mathsf{X})$ is a $C^*$-algebra
with the pointwise operations and the sup-norm.
It is well-known that the Fourier transform is an injective
*-homomorphism of $L^1(G)$ into $\mathfrak{C}_0(\mathsf{X})$
whose range contains the subalgebra $\mathfrak{C}_{00}(\mathsf{X})$
of all cofinite elements. Similar statement holds true for
any of the convolution algebras $C(G), L^p(G)$ $(1\leq p\leq \infty)$.
For a $\sigma_0 \in \mathsf{X}$ let
\[
\mathfrak{I}_{\sigma_0} =\{ (E_\sigma)_{\sigma\in \mathsf{X}} \, :
                         \, E_\sigma =0 \, \, (\sigma\neq \sigma_0)\}.
\]
Under the Fourier transform, the minimal ideals of any of the previously
mentioned convolution algebras are in a one-to-one correspondence with
the ideals $\mathfrak{I}_\sigma$ $(\sigma\in \mathsf{X})$.
Therefore, the minimal ideals are isomorphic to full matrix algebras and
they are algebraically orthogonal to each
other (i.e. the product of any two of them is $\{0\}$).
The structure theory of group algebras tells us that in
$C(G), L^p(G)$ $(1\leq p< \infty)$ the algebraic direct sum
of the minimal ideals is dense in the corresponding norm topology
(see \cite[(28.39) Theorem and (27.39) Theorem]{HeRo}).

\begin{proof}[Proof of Theorem]
Let us suppose first that $p\neq 2$. In this case the form of isometric
automorphisms of $L^p(G)$ is well-known. For every such automorphism
$\phi$ there exist a continuous group character $\tau :G \to \mathbb
T$ and a (bi)continuous group automorphism
$\varphi: G\to G$ such that
\begin{equation}\label{E:formauto}
\phi(f)(x)=\tau(x)f(\varphi(x)) \qquad (x\in G, f\in L^p(G))
\end{equation}
(see \cite[Theorems 2,3]{Stri} and note that by the uniqueness of the
Haar measure every continuous automorphism of $G$ is measure
preserving). In particular, it follows that
$\phi$ maps $C(G)$ onto $C(G)$ and $\phi$ is a surjective isometry with
respect to the sup-norm.

Let $\psi:L^p(G) \to L^p(G)$ be a local
isometric automorphism of $L^p(G)$.
By \eqref{E:formauto}, $\psi_{|C(G)}$ is a local surjective isometry of
$C(G)$ and our Lemma yields that it is a surjective
isometry. Using Banach-Stone theorem, we have
a continuous function $t: G\to \mathbb{T}$ and a homeomorphism
$g: G\to G$ such that
\begin{equation}\label{E:cgform}
\psi(f)(x)=t(x) f(g(x)) \qquad (x\in G, f\in C(G)).
\end{equation}
Considering the local form of $\psi$ at the function $f=1$, we obtain
that
$t$ is a character. Pick different points $x,y \in G\setminus \{1\}$.
Then $x,y, xy$ are pairwise different and so are $g(x), g(y), g(xy)$.
It is not hard to construct a nonnegative continuous function $f$ on $G$
with the property that $f^{-1}(0)=\{g(x),g(y)\}$, $f^{-1}(1)=\{g(xy)\}$.
By the local form of $\psi$ it follows that there exist a continuous
character $\tau :G \to \mathbb{T}$ and a continuous group
automorphism $\varphi:G \to G$ such that
\begin{equation}\label{E:osszehas}
t \cdot f\circ g=\tau \cdot f\circ \varphi.
\end{equation}
Taking absolute value we arrive at $f\circ g=f\circ \varphi$ and then
we deduce $\{ g(x), g(y)\}=\{ \varphi(x), \varphi(y) \}$ and
$g(xy)=\varphi(xy)$. Since $\varphi$ is an automorphism, we have
either $g(xy)=g(x)g(y)$ or $g(xy)=g(y)g(x)$. Similarly, we can
prove that $g(1)=1$ and $g(x^2)=g(x)^2$. Therefore, the homeomorphism
$g: G\to G$ has the property that for every $x,y\in G$ we have
either $g(xy)=g(x)g(y)$ or $g(xy)=g(y)g(x)$. By
\cite[Theorem 2]{Sco} it follows that $g$ is either an automorphism
or an antiautomorphism of $G$.
Apparently, this implies that $\psi$ is either an automorphism or an
antiautomorphism of $C(G)$.

Let $p<\infty$. Since in this case $C(G)$ is dense in $L^p(G)$ and $g$,
being a continuous automorphism or antiautomorphism of $G$, is measure
preserving, it follows easily that the formula \eqref{E:cgform} holds
true also for the elements of $L^p(G)$.
If $p=\infty$, then we can arrive at the same conclusion by using the
fact that $\psi$ is an isometry with respect to any $p$-norm and
$L^\infty (G) \subset L^p(G)$.

Let us stop here for a while and deal with the case $p=2$.
In this case we do not
have the form \eqref{E:formauto} of isometric automorphisms but we can
use Plancherel's theorem instead. It says that via the Fourier
transform, $L^2(G)$ is isometric and isomorphic to
\begin{equation}\label{E:C2}
\mathfrak{C}_2(\mathsf{X})=\{
(A_\sigma)_{\sigma\in \mathsf{X}} \, : \,
(\sum_{\sigma\in \mathsf{X}} d_{\sigma} \| A_\sigma
\|_2^2)^{1/2} <\infty \}.
\end{equation}
Here, $\| .\|_2$ denotes the
Hilbert-Schmidt norm of matrices and on $\mathfrak{C}_2(\mathsf{X})$ we
consider the norm suggested in \eqref{E:C2}. For what remains
we need to know how the isometric Jordan automorphisms of the
$H^*$-algebra $\mathfrak{C}_2(\mathsf{X})$ looks like.
Recall that a linear map $\phi$
between algebras $\mathcal A$ and $\mathcal B$ is called a Jordan
homomorphism if it satisfies
\[
\phi(x^2)=\phi(x)^2 \qquad (x\in \mathcal A),
\]
or equivalently
\[
\phi(xy+yx)=\phi(x)\phi(y)+\phi(y)\phi(x) \qquad (x,y\in \mathcal A).
\]
From the proof of \cite[Theorem]{MolPAMS} it follows that for any
isometric Jordan automorphism $\phi$ of $\mathfrak{C}_2(\mathsf{X})$
there
are a bijection $\alpha :\mathsf{X} \to \mathsf{X}$ and unitary matrices
$U_\sigma \in M_{d_\sigma} (\mathbb C)$ such that $\phi$ is of the form
\begin{equation}\label{E:formC2}
\phi( (A_\sigma)_\sigma)=(U_{\alpha(\sigma)} A_{\alpha(\sigma)}^{[T]}
U_{\alpha(\sigma)}^*)_\sigma
\end{equation}
Here, ${}^{[T]}$ denotes that on the corresponding "coordinate" space
we might have to take transpose.
Now, let $\psi$ be a local isometric automorphism of
$\mathfrak{C}_2(\mathsf{X})$. Clearly,
$\psi$ sends idempotents to idempotents. Like on matrix algebras,
it is now a standard argument to show that $\psi$ is a
Jordan homomorphism on the subalgebra
$\mathfrak{C}_{00}(\mathsf{X})$ of $\mathfrak{C}_2(\mathsf{X})$
consisting of all cofinite elements (see, e.g. \cite[Theorem
2]{MolStud}).
By the continuity of $\psi$ and the density of
$\mathfrak{C}_{00}(\mathsf{X})$ in $\mathfrak{C}_2(\mathsf{X})$ we
infer that $\psi$ is a Jordan homomorphism. Let $\sigma_0 \in
\mathsf{X}$.
The preimage of the minimal ideal $\mathfrak{I}_{\sigma_0}$ under $\psi$
is a finite dimensional Jordan ideal. Suppose that
$\psi^{-1}(\mathfrak{I}_{\sigma_0})\neq \{ 0\}$. Since the full
matrix algebras
are simple Jordan algebras, it follows easily that this preimage is of
the form \begin{equation}\label{E:preim}
\psi^{-1} (\mathfrak{I}_{\sigma_0})=
                \{ (A_\sigma)_{\sigma\in \mathsf{X}} \, :
                \, A_\sigma =0 \, \, (\sigma\notin \{ \sigma_1, \ldots,
                \sigma_n\}) \}.
\end{equation}
By the local form of $\psi$ we get that
in \eqref{E:preim} we have $d_{\sigma_1}=\ldots
=d_{\sigma_n}=d_{\sigma_0}$ and $n=1$.
We obtain $\psi(\psi^{-1} (\mathfrak{I}_{\sigma_0}))=
\mathfrak{I}_{\sigma_0}$. To sum up, if a minimal ideal contains a
nonzero element of the range of $\psi$, then this minimal ideal is
included in the range. Now,
the result \cite[(28.2) Theorem]{HeRo} says that if $G$
is infinite, then the cardinality of $\mathsf{X}$ is equal to the
topological weight of $G$. Since $G$ is metrizable and compact, it
follows that it is second countable. We have
$\mathsf{X}=\{ \sigma_n\}_n$. For every $n$, let $x_n\in
\mathfrak{I}_{\sigma_n}$ be its unit. Consider the element
\[
x=\sum_n\frac{1}{n \sqrt{d_{\sigma_n}}}x_n
\]
and pick a minimal ideal, say $\mathfrak{I}_{\sigma_1}$.
By the local form of
$\psi$ we clearly have $\psi(x)\mathfrak{I}_{\sigma_1}\neq \{0\}$.
This implies that there is an index $n$ such that
$\psi(x_n)\mathfrak{I}_{\sigma_1}\neq \{0\}$. But by the local form of
$\psi$ once again, it follows that $\psi(x_n)$ is the unit of some
minimal ideal. This result in $\psi(x_n)\in
\mathfrak{I}_{\sigma_1}\neq \{0\}$. By what we have proved
before, we conclude that
$\mathfrak{I}_{\sigma_1}$ is included in the range of $\psi$. Since this
minimal ideal was arbitrary and the range of $\psi$ is closed, we obtain
that $\psi$ is surjective and hence, it is an isometric Jordan
automorphism of $\mathfrak{C}_2(\mathsf{X})$.

Now, we are in a position to complete the proof of our theorem.
Let $p\neq 2$. Since $\psi$ is an isometric automorphism or
antiautomorphism of $L^p(G)$, we are done in the commutative case.
If $G$ is finite, then $G$ is discrete and hence there is an
injective nonnegative continuous function on $G$. If we put this
function
into \eqref{E:osszehas}, we obtain that $g$ is an automorphism of $G$
which implies that $\psi$ is an isometric automorphism of $L^p(G)$.
Next suppose that $d_{\sigma_n} \geq 3$ for some $\sigma_n \in
\mathsf{X}$. If $\psi$ is an antiautomorphism, then by its
form given in \eqref{E:cgform}, $\psi$ can be extended from $C(G)$
to an isometric antiautomorphism of $L^2(G)$. By the local
property of $\psi$ we obtain that at the elements of $C(G)$, $\psi$
agrees
with some isometric automorphism of $L^p(G)$ which, keeping its form
in mind, can be extended or restricted to an isometric automorphism of
$L^2(G)$.
The form of isometric automorphisms and antiautomorphisms of $L^2(G)$
can be easily obtained from \eqref{E:formC2}. In fact, in the former
case there is no transpose in \eqref{E:formC2} at all, while in the
latter one we must take
transpose on every "coordinate" space. Therefore, assuming that $\psi$
is an antiautomorphism we would obtain that every matrix
$A\in M_{d_{\sigma_n}}(\mathbb C)$ is unitarily equivalent to its
transpose.
But this is a contradiction. Indeed, one can verify quite easily that
the matrix $A=(a_{ij})$, where $a_{i+1,i}=i$ $(i=1,\ldots ,
d_{\sigma_0}-1)$ and the other entries are all zero, is not unitarily
equivalent to its transpose.
Finally, in the case when $p=2$, the proof can be completed in the same
way observing that in one-dimension the transpose does not matter.
\end{proof}

\begin{remark}
We conclude the paper with some comments on our result.

First, consider the topological condition in Theorem, i.e. the first
countability of $G$.
We suspect that it is essential in the case $p\neq 2$
as well but as for $p=2$, we have the following counterexample.
If $G$ is a compact group with weight greater than $\aleph_0$,
then by \cite[(28.2) Theorem]{HeRo} we have a dimension $d$ which
appears
uncountably many times among the $d_\sigma$-s. Without serious loss of
generality we can restrict our attention to the corresponding part of
$\mathfrak{C}_2(\mathsf{X})$. So, let us consider the $H^*$-algebra
\[
\{ f:\Lambda \to M_d(\mathbb C) \, : \,
      (\sum_{\lambda\in \Lambda} \| f(\lambda)\|_2^2)^{1/2} <\infty \}
\]
of all Hilbert-Schmidt functions, where
$\Lambda$ is an uncountable set.
Let $\Lambda'$ be a proper subset of $\Lambda$ with a bijection $\alpha
:\Lambda' \to \Lambda$. Define $\psi(f)$ by
\[
\psi(f)(\lambda)=
\cases   f(\alpha(\lambda)) \quad \text{if} \quad \lambda\in \Lambda' \\
         0 \qquad \quad \, \, \, \, \, \text{if} \quad \lambda \in
         \Lambda \setminus \Lambda' \endcases.
\]
If $f$ is a function from our collection, then it takes
nonzero values only on a countable set $\{ \lambda_n\}_n$. Let $\beta$ be a
bijection from $\Lambda \setminus \{ \alpha^{-1}(\lambda_n)\}_n$ onto
$\Lambda \setminus \{ \lambda_n\}_n$. Define a bijection $\gamma
:\Lambda \to \Lambda$ by
\[
\gamma(\lambda)=
\cases   \alpha(\lambda) \quad \text{if} \quad \lambda \in \{
         \alpha^{-1}(\lambda_n)\}_n \\
         \beta(\lambda) \quad  \text{if} \quad \lambda \in \Lambda
         \setminus \{ \alpha^{-1}(\lambda_n)\}_n\endcases.
\]
It is easy to see that $\psi(f)=f\circ \gamma=:\phi(f)$, and
$\phi$ is an isometric isomorphism of the algebra of all Hilbert-Schmidt
functions. Therefore, $\psi$ is a nonsurjective local isometric
automorphism.

Our second remark is the following. In the proof of our theorem we have
seen that if $p\neq 2$, every local isometric automorphism $\psi$ of
$L^p(G)$ is of the form
\begin{equation}\label{E:kin}
\psi(f)=t \cdot f\circ g \qquad (f\in L^p(G))
\end{equation}
where $t:G\to \mathbb T$ is a continuous character and $g:G \to G$ is
either a continuous automorphism or a continuous
antiautomorphism. This implies that $\psi$ is either an automorphism or
an antiautomorphism of $L^p(G)$ regardless the possible additional
properties of $G$. It is a natural question whether the conditions
listed in (i) cover every possibility. We mean that if a compact
metric group is noncommutative, infinite and has only one- and
two-dimensional irreducible representations, does it follow that there
is
a local isometric automorphism of $L^p(G)$ which is not an automorphism.
Unfortunately (or fortunately), the answer to this question is negative.
To see this, pick, for
example, an arbitrary noncommutative finite group $H$ and consider the
product $G=\mathbb T \times H$. If $H$ has only one- and
two-dimensional irreducible representations (e.g. $H$ is the symmetric
group
$\mathfrak{S}_3$), then the same holds true for $\mathbb T\times H$ as
well. Let $\psi$ be a local
isometric automorphism of $L^p(G)$ which is necessarily of the form
\eqref{E:kin}. If $\varphi :\mathbb T \times H \to \mathbb
T \times H$ is a continuous automorphism, then an easy argument using
the connectedness of $\mathbb T$ and the fact that none of its closed
subspaces is homeomorphic to $\mathbb T$ shows that $\varphi$ can be
written in the form
\begin{equation}\label{E:simbum}
\varphi(z, x)=(\alpha(z,x), \beta(x)) \qquad (z\in \mathbb T, x\in H),
\end{equation}
where $\beta $ is an automorphism of $H$. Clearly, we have a similar
form for the antiautomorphims of $H$.
Now, by \eqref{E:osszehas} we know that for every $f\in L^p(G)$ there
exist a continuous character $\tau:G \to \mathbb T$ and a
continuous automorphism $\varphi:G \to G$ such that
\begin{equation}\label{E:lokszi}
t \cdot f\circ g=\tau \cdot f\circ \varphi.
\end{equation}
Let $f$ be the function defined by $f(z,x_i)=i$, where $z\in \mathbb T$
is arbitrary and $x_i$ is the $i$th element of $H$. Considering suitable
$\tau$ and $\varphi$ and comparing the two sides of
\eqref{E:lokszi} we deduce that the second coordinate function in the
decomposition of $g$, which can be given similarly to \eqref{E:simbum},
is an automorphism of $G$.
Since $H$ is noncommutative, we obtain that from the two possibilities,
that $g$ is either an automorphism or an antiautomorphism, the
former one is true. Therefore, we have proved that on $L^p(\mathbb
T\times
H)$ every local isometric automorphism is an automorphism. So there is a
small gap in the $p\neq 2$ part of our theorem which is
not the case with $p=2$ as we shall see immediately. Our last remark
concerning the case $p\neq 2$ is that we would clearly obtain the
conclusion of the theorem if there were an injective nonnegative
function in $C(G)$. More precisely, we would be done even if for every
subset of $G$ with three elements we could guarantee
the existence of a nonnegative function which is injective when
restricted to this subset. Unfortunately, this is not the case even with
the most nicest groups like $\mathbb T$.

Finally,
as for the case $p=2$, we show that the condition that $G$ does not have
any two-dimensional irreducible representation is not only necessary but
also sufficient to reach the desired conclusion. To verify this,
we note that every
$2\times 2$ complex matrix is unitarily equivalent to its transpose
as it was proved in \cite[Remark]{BaMo}. Now, if we suppose that at
least
one of the $d_\sigma$-s, say $d_{\sigma_0}$ is 2, then considering the
map $\psi$ on $\mathfrak{C}_2(\mathsf{X})$
which acts as the identity on the
factors $M_{d_\sigma}(\mathbb C)$, $\sigma \neq \sigma_0$ and acts as
the transposition on $M_{d_{\sigma_0}}(\mathbb C)$, we have a local
isometric automorphism of $L^2(G)$ which is not an automorphism.
\end{remark}


\end{document}